\documentclass[12pt]{amsart}
\usepackage{amssymb, latexsym, amsthm}
\newtheorem{Thm}{Theorem}

\newtheorem{Cor}[Thm]{Corollary}

\newtheorem{Lem}[Thm]{Lemma}
\newtheorem{Prop}[Thm]{Proposition}

\newtheorem{Rem}[Thm]{Remark}

\newtheorem{Ques}[Thm]{Question}

\def\st{{such that}}

\newcommand\JaL{Ja\-cquet-Lang\-lands}
\newcommand\JLc{\JaL\ cor\-res\-pon\-dence}

\newcommand\adeles{{ad\`{e}les}}
\def\subrep{sub-re\-pre\-sen\-ta\-tion}
\def\rep{re\-pre\-sen\-ta\-tion}

\def\ie{{i.e.}}
\def\eg{{e.g.}}
\def\cf{{\it cf. \/}}

\def\uni{{\varpi}}

\def\s{\sigma}

\def\A{{\mathbb A}}
\def\C{{\mathbb C}}
\def\F{{\mathbb F}}
\def\HH{{\mathbb H}}
\def\N{{\mathbb N}}
\def\Q{{\mathbb {Q}}}
\def\R{{\mathbb {R}}}
\def\CS{{\mathcal{S}}}
\def\Z{{\mathbb Z}}

\def\B{{\mathcal B}}

\def\GG{{\mathcal G}}

\def\O{{\mathcal{O}}}

\def\Vals{{\mathcal{V}}}

\DeclareMathAlphabet{\mathscr}{OT1}{pzc}{m}{it}

\newcommand\tensor[1][{}]{{\otimes_{#1}}}

\def\ra{{\rightarrow}}
\def\lra{{\longrightarrow}}

\def\minusset{{-}}

\def\sub{\subseteq}
\def\sup{\supseteq}
\def\({\left(}
\def\){\right)}
\def\isom{{\;\cong\;}}

\def\normali{{\lhd}}
\def\co{{\,{:}\,}}
\def\divides{{\,|\,}}
\newcommand\suchthat{{\,:\ \,}}
\newcommand\subjectto{{\,|\ }}

\newcommand\comp[1]{{{#1}^{\operatorname{c}}}}

\DeclareMathOperator{\spec}{spec}

\DeclareMathOperator{\Ker}{Ker}

\DeclareMathOperator{\Aut}{Aut}

\newcommand\JL{{\operatorname{JL}}}
\newcommand{\Sp}{\operatorname{Sp}} %
\DeclareMathOperator{\Cp}{C}

\DeclareMathOperator{\mychar}{char} \DeclareMathOperator{\Cent}{Z}

\newcommand{\Norm}[1][]{{\operatorname{N}_{#1}}}
\DeclareMathOperator{\Gal}{Gal}

\DeclareMathOperator{\Br}{Br}

\newcommand\cond[2][!]{{\operatorname{cond}_{\if!#1\relax\else{\comp{#1}}\fi}(#2)}}

\newcommand\op[1]{{#1^{\operatorname{op}}}}
\newcommand\trans[1]{{{#1}^{\operatorname{t}}}}
\renewcommand\L[2]{{\operatorname{L}^{#1}(#2)}}
\newcommand\M[1][d]{{\operatorname{M}_{#1}}}
\newcommand\GL[1][d]{{\operatorname{GL}_{#1}}}
\newcommand\PGL[1][d]{{\operatorname{PGL}_{#1}}}
\newcommand\PGU[1][d]{{\operatorname{PU}_{#1}}}
\newcommand\PGO[1][d]{{\operatorname{PO}_{#1}}}

\newcommand\PU[1][d]{{\operatorname{PU}_{#1}}}
\newcommand\SL[1][d]{{\operatorname{SL}_{#1}}}
\newcommand\SO[1][d]{{\operatorname{SO}_{#1}}}
\newcommand\SU[1][d]{{\operatorname{SU}_{#1}}}
\newcommand\PSL[1][d]{{\operatorname{PSL}_{#1}}}

\newcommand\abs[2][F]{|{#2}|_{#1}}

\newcommand{\set}[1]{{\{#1\}}}
\newcommand{\card}[1]{{\left|{#1}\right|}}

\newcommand\dimcol[2]{{[{#1}\!:\!{#2}]}}

\newcommand\mul[1]{{#1^{\times}}}
\newcommand\md[2][d]{{\mul{#2}/{\mul{#2}}^{#1}}}

\newcommand\valF{{\nu_0}}

\newcommand\Tref[1]{{Theorem \ref{#1}}}

\newcommand\Pref[1]{{Proposition \ref{#1}}}
\newcommand\Lref[1]{{Lemma \ref{#1}}}
\newcommand\Cref[1]{{Corollary \ref{#1}}}

\newcommand\eq[1]{{(\ref{#1})}}

\newcommand\Eq[1]{{Equation \eq{#1}}}

\newcommand\defin[1]{{\it{#1}}}
\long\def\half#1\halved{{\footnotesize{#1}}}
\long\def\forget#1\forgotten{}

\newcommand\dc[3]{{{#1}\backslash{#2}/{#3}}}
\newcommand\dom[2]{{{#1}\backslash{#2}}}
\newcommand\firsteqnarray[2][0.55cm]{{\par\vspace{0.4cm}\lefteqn{#2}\nopagebreak\vspace{-#1}}}

\newcommand\paper[6]{{{#1},\ {\it{#2}},\ {#3}\ {#4},\ {#5},\ ({#6}).}}
\newcommand\book[4]{{{#1},\ {{#2}},\ {#3},\ {#4}.}}

\newcommand\absdot[1][]{\abs[#1]{\;\!\cdot\:\!}}
\newcommand\dd[4][!]{\abs[#2]{\:\!\det{\if!#1\relax\:\!\!\else(#1)\fi}}^{#3} #4}

\newcommand\automrep[1][G']{{\L2{\dom{{#1}(k)}{{#1}(\A)}}}}

\newcommand\con[3]{{{#1}({#2},{#3})}}

\newcommand{\oper}[1]{{\operatorname{#1}}}
\def\GAT{{G_0}}

\newcommand\ver[2][Changed:]{{}}

\newif\ifXY
\XYfalse

\ifXY
\usepackage{xy}
\fi
\ifXY
\xyoption{all}
\fi

\begin{document}

\title
[Division Algebras and Isospectrality]
{Division Algebras and Non-Commensurable Isospectral Manifolds}

\def\HUJI{Inst. of Math., Hebrew Univ., Givat Ram, Jerusalem 91904,
Israel}
\def\BIU{Dept. of Math., Bar-Ilan University, Ramat-Gan 52900,
Israel}
\def\YALE{Dept. of Math., Yale University, 10 Hillhouse Av., New-Haven CT
06520, USA}

\author{Alexander Lubotzky}
\address{\HUJI}
\email{alexlub@math.huji.ac.il}

\author{Beth Samuels}
\address{\YALE
}
\email{beth.samuels@yale.edu}

\author{Uzi Vishne }
\address{\BIU}

\email{vishne@math.biu.ac.il}

\renewcommand{\subjclassname}{
      \textup{2000} Mathematics Subject Classification}

\date{Submitted: Dec. 30, 2004}

\maketitle

\begin{abstract}
A.~Reid \cite{Reid} showed that if $\Gamma_1$ and $\Gamma_2$ are
arithmetic lattices in $G = \PGL[2](\R)$ or in $\PGL[2](\C)$ which
give rise to isospectral manifolds, then $\Gamma_1$ and $\Gamma_2$
are commensurable (after conjugation).
We show that for $d \geq 3$ and $\CS = \PGL[d](\R) / \PGO[d](\R)$,
or $\CS = \PGL[d](\C) / \PU[d](\C)$, the situation is quite
different: there are arbitrarily large finite families of
isospectral non-commensurable compact manifolds covered by $\CS$.

The constructions are based on the arithmetic groups obtained from
division algebras with the same ramification points but different
invariants.
\end{abstract}

\section{Introduction}

Let $X$ be a compact Riemannian manifold
and $\Delta = \Delta_X$ its Laplacian. The spectrum of $X$,
$\spec(X)$, is the multiset of eigenvalues of $\Delta$ acting on
$\L2{X}$.
This is a discrete subset of $\R$. If $Y$ is
another such object, we say that $X$ and $Y$ are
\defin{isospectral} (or, sometimes, cospectral) if $\spec(X) =
\spec(Y)$.
The problem of finding isospectral non-isomorphic objects has a
long history, starting with Marc Kac's seminal paper \cite{Kac},
``Can one
hear the shape of a drum?'' (see \cite{Gordon} and the references
therein).

The most powerful and general method of obtaining such pairs is
due to Sunada \cite{Sunada}: he showed that if $X_0$ covers $X$
with a finite Galois group $H$, \ie\ $X = X_0/H$, and if $H_1$ and
$H_2$ are subgroups of $H$ such that for every conjugacy class $C
\subset H$, $\card{C \cap H_1} = \card{C \cap H_2}$, then $X_1 =
X_0/H_1$ and $X_2 = X_0/H_2$ are isospectral. Of course, one still
has to determine whether $X_1$ and $X_2$ are isomorphic (for
example, $X_1$ and $X_2$ will be isomorphic if $H_1$ is conjugate
to $H_2$, but not only in such a case). Sunada's method has been
implemented in many situations (\eg\ \cite{Br1}, \cite{Br2}), and,
in particular, has led to the solution of the original problem of
M.~Kac (see \cite{GWW}).

A common feature to all the isospectral pairs, $X_1$ and $X_2$,
constructed by Sunada's method, is that they are commensurable,
that is, $X_1$ and $X_2$ have a common finite cover $X_0$. So,
while $X_1$ and $X_2$ are not isomorphic, they are `virtually
isomorphic'.

Sunada implemented his method
to obtain non-isomorphic isospectral Riemann surfaces.  He gave
examples of torsion free cocompact lattices $\Gamma \normali
\Gamma_0$ in $G = \PGL[2](\R)$ with $\Gamma_0 / \Gamma \isom H$
and $H_1 = \Gamma_1/ \Gamma$, and $H_2 = \Gamma_2 / \Gamma$, as
above. His method implies that the $G$-\rep{s}
$\L2{\dom{\Gamma_1}{G}}$ and $\L2{\dom{\Gamma_2}{G}}$ are
isomorphic, which demonstrates that $\dom{\Gamma_1}{\HH^2}$ is
isospectral to $\dom{\Gamma_2}{\HH^2}$. Here $\HH^2 =
\PGL[2](\R)/K$ is the associated symmetric space, and $K$ is the
maximal compact subgroup $\PGO[2](\R)$ of $G$.

Vigneras was the first to present examples of isospectral Riemann
surfaces in \cite{Vig}, where she used the theory of quaternion
algebras.  Her examples are also commensurable to each other. In
fact, it is still an open problem whether every two isospectral
Riemann surfaces are commensurable, or equivalently whether every
two cocompact lattices $\Gamma_1$ and $\Gamma_2$ in $G =
\PGL[2](\R)$, such that $\dc{\Gamma_1}{G}{K}$ and
$\dc{\Gamma_2}{G}{K}$ are isospectral, are commensurable (after
conjugation, \ie\ for some $g \in G$,
$\dimcol{\Gamma_2}{g^{-1}\Gamma_1 g \cap \Gamma_2} < \infty$).
But, Alan Reid \cite{Reid} showed that if $\Gamma_1$ and
$\Gamma_2$ are such
{\emph {arithmetic}} lattices, then
indeed $\Gamma_1$ and $\Gamma_2$ are commensurable. He also showed
a similar result for $G = \PGL[2](\C)$.

Our main result shows that the situation is quite different for
$G = \PGL[d](\R)$ or $\PGL[d](\C)$ if $d \geq 3$.

\begin{Thm}\label{isomanifold}
Let $F= \R$ or $\C$, $G = \PGL[d](F)$ where $d \geq 3$, $K\leq G$
a maximal compact
subgroup,
and $\CS = G/K$ the associated symmetric space. Then for every $m
\in \N$, there exists a family of $m$ torsion free cocompact
arithmetic lattices $\set{\Gamma_i}_{i=1,\dots,m}$ in $G$, such
that
$\dom{\Gamma_i}{\CS}$ are isospectral and not commensurable.
\end{Thm}

Taking $K = \PGO[d](\R)$ in the real case and $K =
\PGU[d](\C)$ in the complex case, the quotients
$\dc{\Gamma_i}{G}{K}$ are isospectral compact manifolds covered by
the symmetric space $\CS = \PGL[d](F)/K$ which is equal to
$\SL[d](\R)/\SO[d](\R)$ or $\SL[d](\C) / \SU[d](\C)$,
respectively. The covering map $\CS \ra \dom{\Gamma_i}{\CS}$ is a
local isomorphism, since the $\Gamma_i$ are torsion free.

This result should be compared with Spatzier's construction
\cite{Spatz} of
isospectral locally symmetric spaces of high rank.
However, his isospectral examples are always commensurable to each
other.

\medskip

Let us now outline the method of proof. The following is well
known (see for example \cite{Pesce}).
\begin{Prop}\label{isospecdef}
Let $G$ be a semisimple group,
$K \leq G$ a maximal
compact subgroup, and $\Gamma_1,\Gamma_2 \leq G$ discrete
cocompact subgroups. If
\begin{equation}\label{eqgeneral}
\L2{\dom{\Gamma_1}{G}} \isom
\L2{\dom{\Gamma_2}{G}}\end{equation}
as (right) $G$-representations,
then $\dc{\Gamma_1}{G}{K}$ and $\dc{\Gamma_2}{G}{K}$ are
isospectral quotients of $G/K$ (indeed, they are even strongly
isospectral, \ie\ with respect to the higher dimensional
Laplacians).
\end{Prop}

In some special cases, the inverse of \Pref{isospecdef} is known
to be true (\cf\ \cite{Pesce} and the references therein, and
\cite{Demir}).

To prove Theorem \ref{isomanifold} for $G = \PGL[d](F)$ ($F = \R$
or $\C$), we choose the discrete subgroups $\Gamma_i$ to be
arithmetic lattices of inner forms, as follows.  Let $k$ be a
global field and $\set{k_{\nu}}$ the completions with respect to
its valuations. Recall that $\Br(k_{\nu}) \isom \Q/\Z$ for a
non-archimedean valuation, while $\Br(\R)\isom \frac{1}{2}\Z/\Z$
and $\Br(\C) = 0$, where $\Br()$ denotes the Brauer group of a
field. For a division algebra $D$, $[D] \in \Br(k)$, the value
associated to $D \tensor[k] k_{\nu}$ is called the
`local invariant' at ${\nu}$.
By Albert-Brauer-Hasse-Noether theorem,
$[D] \mapsto ([D \tensor[k] k_{\nu}])_{\nu}$ defines an injective
map for the Brauer groups
\begin{equation}\label{Breq}
\Br(k) \lra \oplus  \Br(k_{\nu}), \end{equation} and the image of
this map is composed of vectors whose sum is zero. Over local (and
thus also global) fields, the degree of a division algebra is
equal to its exponent (which is the order of its equivalence class
in the Brauer group).

Fix  $d > 1$
and let $T = \set{\theta_1,\dots,\theta_t}$ be a finite non-empty
set of non-archimedean valuations of $k$. Let
$a_1,\dots,a_t,b_1,\dots,b_t \in \N$ be \st\ for every $j =
1,\dots,t$, $0 \leq a_j,b_j < d$, $(a_j,d) = (b_j,d) = 1$, and
 $\sum{a_j} \equiv \sum{b_j} \equiv 0 \pmod{d}$. Let $D_1$
(resp. $D_2$) be the unique division algebra of degree $d$ (\ie\
dimension $d^2$) over $k$ which ramifies exactly at $T$, and whose
invariant at
$k_{\theta_j}$ is $a_j/d$ (resp. $b_j/d$).
Thus $D_1$ and $D_2$ split at every $\nu \not \in T$. Let $G'_i$
be the algebraic group $\mul{D_i}/\mul{Z}$ ($i = 1,2$) where $Z$
is the center, and $G = \PGL[d]$. Since $G'_i(k_{\nu})$ splits for
every $\nu \not \in T$, the groups $G'_1(\A_{\comp{T}})$ and
$G'_2(\A_{\comp{T}})$ are equal where $\A_{\comp{T}} =
\prod'_{\nu\not\in T} k_{\nu}$.  We identify these groups with
$\GAT = G(\A_{\comp{T}})$.  Here, $\A$ is the ring of \adeles\ and
$\prod'$ denotes the restricted product, namely the vectors
$(x_{\nu})$ such that $\nu(x_{\nu}) \geq 0$ for almost every
$\nu$.

A basic result is
\begin{Thm}\label{basic}
Let $\GAT = \PGL[d](\A_{\comp{T}})$ be the group defined above,
and let $\Delta_i =
G'_i(k)$ ($i = 1,2$), cocompact
lattices in $\GAT$. The spaces
$\L2{\dom{\Delta_1}{\GAT}}$ and $\L2{\dom{\Delta_2}{\GAT}}$ are
isomorphic as $\GAT$-\rep{s}.
\end{Thm}

The result follows by comparing trace formulas.  This was done
explicitly for the case where $d$ is prime in
\cite[Theorem~1.12]{GJ}
and in \cite[Theorem~54]{Bump2}. It seems that the general case is
also known to experts, but we were unable to locate a reference.
For the sake of completeness we show in Section \ref{later} how
the result can be deduced (for arbitrary $d$) from the global
Jacquet-Langlands correspondence, as described in \cite{HT}.
(The proof there is given for the characteristic zero case, and we
refer the reader to \cite[Remark~1.6]{paperI} for some details
about the positive characteristic case).

We should stress that in the current paper, we are using only the
characteristic zero case. We prefer to give here the more general
formulation, preparing for a subsequent paper \cite{paperIV}, in
which we use analogous ideas over a local field of positive
characteristic to construct isospectral simplicial complexes and
isospectral Cayley graphs of some finite simple groups.

\ver{4.4} One deduces the isospectrality in Theorem
\ref{isomanifold} by applying strong approximation to
\Tref{basic}. This is done in Section \ref{thm12}. To prove the
non-commensurability, we first show (\Tref{noncommens}) that two
division $k$-algebras $D_1$ and $D_2$ give rise to commensurable
arithmetic lattices in $\PGL[d](F)$ iff $D_1(k)$ and $D_2(k)$ are
isomorphic or anti-isomorphic (as rings, rather than as
$k$-algebras). We then analyze when two such division algebras are
isomorphic (or anti-isomorphic) as rings, and show that we can
produce as many examples as we need of non-isomorphic division
rings.  This is done in Section \ref{sec:nonc}.

Finally, we mention that our work is very much in the spirit of
the paper of Vigneras \cite{Vig}, who used a similar idea to find
isospectral Riemann surfaces. However, as mentioned above, her
examples, which are lattices in $\PGL[2](\R)$, are commensurable,
as they should be by Reid's theorem. The case $d = 2$ differs from
those of $d \geq 3$ in that a global division algebra of degree
$2$ is determined by its ramification points, while for $d \geq 3$
there are non-isomorphic division algebras with the same
ramification points.

We thank G.~Margulis for helpful discussions, A.S.~Rapinchuk for
his help with \Tref{noncommens}, and D.~Goldstein and R.~Guralnick
for simplifying the proof of \Lref{subring}. The authors also
acknowledge support of grants by the NSF and the U.S.-Israel
Binational Science Foundation.

\section{Isospectrality}\label{thm12}

Although for the proof of \Tref{isomanifold} we may assume
characteristic zero, we are stating and proving the results of
this section for positive characteristic as well. Therefore, we
obtain \Tref{isomanifold} together with:
\begin{Thm}\label{isoComplex}
Let $F$ be a local field of positive characteristic, $G =
\PGL[d](F)$ where $d \geq 3$, $K$ a maximal compact subgroup and
$\B_d(F) = G/K$ the associated Bruhat-Tits building. Then for
every $m \in \N$ there exists a family of $m$ torsion free
cocompact arithmetic lattices $\set{\Gamma_i}_{i = 1,\dots,m}$ in
$G$, such that the finite complexes $\dom{\Gamma_i}{\B_d(F)}$
 are isospectral and not commensurable.
\end{Thm}

Notice that Theorems \ref{isomanifold}  and \ref{isoComplex} cover
the archimedean local fields, and the non-archimedean local fields
of positive characteristic. Our methods do not apply to
non-archimedean local fields $F$ of zero characteristic, where
there are no cocompact lattices of inner type in $\PGL[d](F)$ ($d
\geq 3$). There are cocompact lattices of outer type in these
groups; however we leave open the following problem:
\begin{Ques}
Let $F$ be a non-archimedean local field of zero characteristic,
$G = \PGL[d](F)$ where $d \geq 3$, $K$ a maximal compact subgroup
and $\B_d(F) = G/K$ the associated Bruhat-Tits building. Do there
exist torsion free cocompact arithmetic lattices
$\set{\Gamma_i}_{i = 1,2}$ in $G$, such that the finite complexes
$\dom{\Gamma_i}{\B_d(F)}$
 are isospectral and not commensurable?
\end{Ques}

Let $k$ be a global field, and assume it has at most one
archimedean valuation.
Let $\nu_0$ denote the archimedean valuation if $\mychar k = 0$,
and a valuation of degree $1$ if $k$ is a function field (for $k =
\F_q(t)$, $\nu_0$ is either the minus-degree valuation or is
determined by a linear prime). Let $F = k_{\nu_0}$.
(For $F = \R$ or $F = \C$, we take $k = \Q$ or $k = \Q[i]$,
respectively).

Let $\Vals$ denote the set of all valuations of $k$ other than
$\nu_0$. For a valuation $\nu \in \Vals$, $\O_{\nu}$ is the ring
of integers in the completion $k_{\nu}$, and $P_{\nu}$ is the
valuation ideal.

Let $T \sub \Vals$ be a finite subset, and define $D_i$ and $G_i'$
as in the introduction ($i = 1,2$). Let $G$ denote the algebraic
group $\PGL[d]$. Notice that $\nu_0 \not \in T$, so in particular
$G_i'(F) = G(F)$ ($i = 1,2$). Let $\A_T$ denote the direct product
$\prod_{\theta \in T} k_{\theta}$, so that $\A = \A_T \times
\A_{\comp{T}}$ is the ring of \adeles\ over $k$.

Let  $R_{0} = \set{x \in k \suchthat \nu(x)\geq 0 \quad \mbox{for
every $\nu \neq \valF$}}$. (In the cases $k = \Q$ or $k = \Q[i]$,
with $\nu_0$ as the archimedean valuation, we obtain $R_0 = \Z$ or
$R_0 = \Z[i]$, respectively.)
Then, for every $\theta
\in T$, we can choose
a uniformizer $\uni_{\theta}$
such that $\nu(\uni_{\theta}) = 0$ for every $\nu \not \in \set{
\nu_0,\theta}$ (so in particular, $\uni_{\theta} \in R_0$). In
this situation,
\begin{equation}\label{R0T}
R_{0,T} = \set{x \in k \suchthat \nu(x)\geq 0 \quad \mbox{for
every $\nu \not\in T \cup \set{\valF}$}}
\end{equation}
is equal to
$R_0[\set{\uni_\theta^{-1}}_{\theta\in T}]$.
Fix $i  \in \set{1,2}$. Since division algebras over global fields
are cyclic, there is a Galois field extensions $k_i/k$ of
dimension $d$ contained in $D_i(k)$.  If $\phi \in \Gal(k_i/k)$ is
a generator, there is an element $b_i \in \mul{k}$ such that $D_i
= k_i[z \subjectto zaz^{-1} = \phi(a),\,z^d = b_i]$. Taking an
integral basis, one finds an order $O_i \sub k_i$ which is a
cyclic extension (of rings) of $R_{0,T}$. The constant $b_i$ can
be chosen to be in the multiplicative group generated by the
uniformizers $\uni_{\theta}$. Thus, $b_i$ is invertible in
$R_{0,T}$.  So $O_i[z]$ (with the relations as above) is an
Azumaya algebra of rank $d$ over the center $R_{0,T}$, and $D_i =
k \tensor[R_{0,T}] O_i[z]$. The group of invertible elements in
$O_i[z]$ (again with the above relations), modulo its center, is
$G'_i(R_{0,T})$.

If an algebraic group $G'$ is defined over a ring $R$ and $0 \neq
I \normali R$ is an ideal, we have the principal congruence
subgroups
\begin{equation}\label{condef}
\con{G'}{R}{I} = \Ker(G'(R) \ra G'(R/I)).
\end{equation}
The congruence subgroups $\con{G'_i}{R_{0,T}}{I}$
are discrete in $G(F)$,
and co-compact if $T \neq \emptyset$ \cite{PR}.

Let $0 \neq I \normali R_{0,T}$.
Define a function $r \co \Vals \minusset (T\cup\set{\valF}) \ra
\N\cup \set{0}$ by setting
\begin{equation}\label{rdef}
r_{\nu} = \min\set{\nu(a)\suchthat a\in I},
\end{equation}
and let
\begin{equation}\label{Urdef}
U^{(r)} = \prod_{\nu \in \Vals - (T\cup \set{\nu_0})} \con{G}{\O_{\nu}}{P_{\nu}^{r_{\nu}}},
\end{equation}
an open compact subgroup of $G(\A_{T^{c}-\set{\nu_0}})$. For
almost every $\nu$, $r_{\nu} = 0$, and then $P_{\nu}^{0} =
\O_{\nu}$ and $\con{G}{\O_{\nu}}{P_{\nu}^{r_{\nu}}} = G(\O_{\nu})$
by definition.
Let
\begin{equation}\label{Uidef}
U_i = U^{(r)} G'_i(\A_{T}),
\end{equation}
an open compact subgroup of $G'_i(\A)$.

\begin{Lem}\label{finqd}
The quotient $G'_i(\A) / ( G_i'(k)G(F)U_i)$ is a $d$-torsion
finite abelian group, which is independent of $i$.
\end{Lem}
\begin{proof}
Recall the reduced norm map $\mul{D_i(k)} \ra \mul{k}$, which
coincides with the determinant over splitting fields of $k$. Let
$\GG_i$ denote the simply connected cover of $G_i'$, \ie\ $\GG_i$
is the subgroup of elements of reduced norm $1$ in $\mul{D_i}$.
For every field $k_1 \sup k$, the image of the covering map
$\Psi\co \GG_i(k_1) \ra G_i'(k_1)$ is co-abelian in $G'_i(k_1)$.
The strong approximation theorem \cite{PR} (see also
\cite[Subsection~3.2]{paperI}) applies to $\GG_i$, and so
$\GG_i(k)\GG_i(F)U = \GG_i(\A)$ where
$$U= \prod_{\nu \in \Vals - (T\cup \set{\nu_0})}
\con{\GG_i}{\O_{\nu}}{P_{\nu}^{r_{\nu}}} \times
\GG_i(\A_{T}),$$
and $\GG_i(\O_{\nu_0}) U$ is open compact in $\GG_i(\A)$.

The reduced norms $G'_i(k_{\nu}) \ra \md[d]{k_{\nu}}$ is onto
whenever $G'_i$ splits at $k_{\nu}$ (as $G'_i(k_{\nu}) =
\PGL[d](k_{\nu})$), and also whenever $k_{\nu}$ is non-archimedean
(since over local fields the norm map in a division algebra is
onto, \cite[Chapter~17] {Pierce}). In the second case, the
restriction to $G'_i(\O_{\nu}) \ra \md[d]{\O_{\nu}}$ is also onto.
Thus, the norm map $G'_i(\A) \ra \md[d]{\A}$ is onto, and induces
an isomorphism $\Phi \co G'_i(\A) / \Psi (\GG_i(\A)) \ra
\md[d]{\A}$.

Since $\Psi$ maps $\GG_i(k)$, $\GG_i(F)$ and $U$ into $G'_i(k)$,
$G(F)$ and $U_i$, respectively, we have $\Psi (\GG_i(\A)) =
\Psi(\GG_i(k)\GG_i(F)U) \sub G'_i(k) G(F) U_i$.
We need to compute the image of $G'_i(k)G(F) U_i$ under $\Phi$.
For a valuation ring $\O$ with maximal ideal $P$ and $r \geq 0$,
let $\O^{(r)}$ denote the subgroup $\set{a \in\mul{\O} \suchthat a
\equiv 1 \pmod{P^{r}}}$. Now, modulo $d$ powers, the determinant
maps $G(F)$ to $\mul{F}$ and
$\con{G}{\O_{\nu}}{P_{\nu}^{r_{\nu}}}$ to $\O_{\nu}^{(r_{\nu})}$.
It follows that $\Phi$ takes $G'_i(k)$ to $\mul{k}\mul{\A}^d$,
$G(F)$ to $\mul{F} \mul{\A}^d$, and $U_i$ to the product
$(\prod_{\nu \in \Vals - T}{\O_{\nu}^{(r_{\nu})}} \times
\prod_{\theta \in T}{\mul{k_{\theta}}})\mul{\A}^d$. We can now
compute that
\begin{eqnarray*}
G'_i(\A) / (G'_i(k)G(F)U_i) & \isom & \mul{\A} /
(\mul{k}\mul{F}\mul{\A}^d \prod_{\nu\in \Vals -(T\cup \set{\nu_0})
}{\O_{\nu}^{(r_{\nu})}}
\prod_{\theta \in T}{\mul{k_{\theta}}} ) \\
    & \isom & \mul{\A_{\comp{T}}} / (\mul{k} \mul{F}\mul{\A_{\comp{T}}}^d
\prod_{\nu\in \Vals -(T\cup \set{\nu_0})}{\O_{\nu}^{(r_{\nu})}} ).
\end{eqnarray*}

Recall that $\mul{k_{\nu}} \isom \Z \times \mul{\O_{\nu}}$, where
the $\Z$ summand is the value group. Dividing the numerator and
denominator in the last expression by $\prod_{\nu \in \Vals -
(T\cup \set{\nu_0})}{\mul{\O_{\nu}}^d}$, the quotient is
\firsteqnarray[0.45cm]{ {\prod^{}_{}}'{\mul{k_{\nu}}} /
   (\mul{k} {\prod_{}}'{\O_{\nu}^{(r_{\nu})}\mul{k_{\nu}}^d})
\isom }
\begin{eqnarray*}
   \qquad
   & \isom &
{\prod^{}_{}}'({\Z \times \mul{\O_{\nu}}}) /
   (\mul{k} {\prod_{}^{}}'({d\Z\times \O_{\nu}^{(r_{\nu})}\mul{\O_{\nu}} })) \\
   \qquad
   & \isom &
{\prod^{}_{}}'({\Z \times \md[d]{\O_{\nu}}}) /
   (\mul{k} {\prod_{}^{}}'({d\Z\times \O_{\nu}^{(r_{\nu})}\md[d]{\O_{\nu}} })) \\
   & \isom &
({\prod^{}_{}}'{\Z } {\prod^{}_{}}{\md[d]{\O_{\nu}}}) /
   (\mul{k} {\prod_{}^{}}'{d\Z} {\prod_{}}{\O_{\nu}^{(r_{\nu})}\md[d]{\O_{\nu}}
    })
\\
   & \isom &
({\prod^{}_{}}'({{\Z/d } \times \mul{\O_{\nu}}/{\O_{\nu}^{(r_{\nu})}\mul{\O_{\nu}}^d}}))
/ \mul{k},
\end{eqnarray*}
where all products are over $\nu \in \Vals - (T\cup \set{\nu_0})$
and $\prod'$ denotes the restricted product (such that
$\nu(a_{\nu}) \geq 0$ almost always).
The quotients
$\mul{\O_{\nu}}/\O_{\nu}^{(r_{\nu})}$ are always finite, and
$r_{\nu} = 0$ for almost all $\nu$, hence
$\mul{\O_{\nu}}/{\O_{\nu}^{(r_{\nu})}\mul{\O_{\nu}}^d} = 0$ for
all but finitely many valuation. In the last expression, $\mul{k}$
embeds in the $\nu$ component by $a \mapsto
(\nu(a),a{\uni_\nu}^{-\nu(a)}{\O_{\nu}^{(r_{\nu})}\mul{\O_{\nu}}^d})$,
and in particular the components with
$\mul{\O_{\nu}}/{\O_{\nu}^{(r_{\nu})}\mul{\O_{\nu}}^d} = 0$ vanish
in the quotient. Therefore we have a
finite product of finite groups,
and the result follows.
\end{proof}

\begin{Prop}\label{basicprop}
Let $0 \neq I \normali R_{0,T}$.
There is an isomorphism
\begin{equation}\label{basicG}
\L2{\dom{\con{G'_1}{R_{0,T}}{I}}{G(F)}} \isom
\L2{\dom{\con{G'_2}{R_{0,T}}{I}}{G(F)}},
\end{equation}
as representations of $G(F) = \PGL[d](F)$.
\end{Prop}
\begin{proof}
We prove this statement by transferring it to its {ad\`{e}{l}ic}
analogue. Let $U_i$ be the group defined in \Eq{Uidef}. Then
$G(\O_{\nu_0}) U_i$ is open compact in $G'_i(\A)$, and
\begin{equation}\label{defGe}
\con{G'_i}{R_{0,T}}{I} = G'_i(k) \cap G(F) U_i, \end{equation} the
intersection taken in $G'_i(\A)$ and then projected to the $G(F)$
component.  Therefore $\dom{\con{G'_i}{R_{0,T}}{I}}{G(F)}$ is
related to $\dc{G'_i(k)}{G'_i(\A)}{U_i}$.

Since $\GAT = G(\A_{\comp{T}})$ by definition, $G'_i(\A)/ U_i =
G'_i(\A_{\comp{T}})/U^{(r)} = \GAT/U^{(r)}$, and
$$\dc{G'_i(k)}{G'_i(\A)}{U_i} \isom \dc{G'_i(k)}{\GAT}{U^{(r)}}.$$
Consequently,
$$\L2{\dom{G_i'(k)}{\GAT}}^{U^{(r)}} \isom
\L2{\dc{G_i'(k)}{\GAT}{U^{(r)}}} \isom
\L2{\dc{G_i'(k)}{G'_i(\A)}{U_i}}.$$

Now let $H_i = G_i'(k)G(F)U_i$. {}From \Eq{defGe} it follows that
$$\dom{\con{G'_i}{R_{0,T}}{I}}{G(F)} \isom
\dc{G_i'(k)}{H_i}{U_i}.$$
By the lemma, $C = G'_i(\A)/H_i$ is finite ($d$-torsion) abelian,
which is independent of $i$.  Since $G(F)$ commutes with $U_i$ and
is contained in $H_i$, we have the decomposition
$$\L2{\dc{G_i'(k)}{G'_i(\A)}{U_i}} \isom \oplus_{\card{C}}
\L2{\dc{G_i'(k)}{H_i}{U_i}}$$
as $G(F)$-\rep{s}. Thus
$$\L2{\dom{G_i'(k)}{\GAT}}^{U^{(r)}} \isom \oplus_{\card{C}}
\L2{\dc{G_i'(k)}{H_i}{U_i}} \isom \oplus_{\card{C}} \L2{\dom{\con{G'_i}{R_{0,T}}{I}}{G(F)}}.$$

By \Tref{basic}, $\L2{\dom{G_1'(k)}{\GAT}} \isom
\L2{\dom{G_2'(k)}{\GAT}}$ as \rep{} spaces, so the result follows.
\end{proof}

Now, let $0 \neq I \normali R_{0,T}$, and take $\Gamma_i =
\con{G'_i}{R_{0,T}}{I}$. By \Pref{isospecdef} and
\Pref{basicprop}, we see that $\dc{\Gamma_{1}}{G(F)}{K}$ and
$\dc{\Gamma_2}{G(F)}{K}$ are isospectral, where $K$ is a fixed
maximal compact subgroup of $G(F)$. If $I$ is small enough, then
$\Gamma_1$ and $\Gamma_2$ are torsion free, and so the projection
$G(F)/K \ra \dc{\Gamma_i}{G(F)}{K}$ is a local isomorphism.

\begin{Rem}
When $T$ is fixed, the choice of $U^{(r)}$ determines the family
of quotients. When ${U^{(r)}}' \sub U^{(r)}$, the quotients
corresponding to $U^{(r)}$ are covered by those corresponding to
${U^{(r)}}'$. In this sense, we find not only families of
isospectral structures, but infinite ``inverse limits'' of such
families.
\end{Rem}

\section{Non-commensurability}\label{sec:nonc}

As in the previous section, we state and prove the results of this
section for arbitrary characteristic, proving Theorems
\ref{isomanifold} and \ref{isoComplex} together.

\ver{4.2}
To show that $\dc{\Gamma_i}{G(F)}{K}$ are not commensurable, we
need the following theorem. Recall that the maximal compact
subgroup $K$ of $\PGL[d](F)$ is taken to be $K = \PGO[d](\R)$ if
$k = \R$, $K = \PGU[d](\C)$ if $k = \C$, and $K = \PGL[d](\O)$
(where $\O$ is the ring of integers of $F = k_{\nu_0}$) if
$\mychar k > 0$.
\ver{4.2} Let $D_1$ and $D_2$ be two central division algebras
over $k$, as in the previous section, and $G'_i$ the corresponding
algebraic groups. Let $T_1$ and $T_2$ be the ramification points
of $D_1$ and $D_2$ (this time not necessarily equal), and let
$R_{0,T_i}$ be the subrings of $k$ defined in \Eq{R0T} (for $T_i$
rather than $T$). Note that if $\s \co D_1(k) \ra D_2(k)$ is an
isomorphism or anti-isomorphism of rings, then $\s$ induces an
automorphism of the center $k$ which acts on the set of
(non-archimedean) valuations of $k$. This maps $T_1$ to $T_2$.
For \Tref{isomanifold} we only need the implication $(1)\implies
(3)$ of the next theorem, however we give the full picture, which
seems to be of independent interest.

\begin{Thm}\label{noncommens}
\ver{4.2} Let $\CS = G(F)/K$ be the building or symmetric space
corresponding to $G(F)$. The following are equivalent:

1. There exists finite index torsion-free subgroups $\Omega_1$ of
$G'_1(R_{0,T_1})$ and $\Omega_2$ of $G'_2(R_{0,T_2})$ such that
the manifolds or complexes $\dom{\Omega_1}{\CS}$ and
$\dom{\Omega_2}{\CS}$ are isomorphic.

2. For every finite index torsion-free subgroup $\Omega_1 \leq
G'_1(R_{0,T_1})$, there exists a finite index torsion-free
subgroup $\Omega_2 \leq G'_2(R_{0,T_2})$ such that
$\dom{\Omega_1}{\CS}$ and $\dom{\Omega_2}{\CS}$ are isomorphic.

\ver{4.2} 3. There is a ring isomorphism or anti-isomorphism $\s
\co D_1(k) \ra D_2(k)$, such that the restriction of $\s$ to the
center $k$ fixes $\valF$ and maps $T_1$ to $T_2$.
\end{Thm}
\begin{proof}
\ver{4.1}
The second assertion trivially implies the first. We
prove $(1) \implies (3)$.  Since $\CS$ is the universal cover of
$\dom{\Omega_i}{\CS}$ and $\Omega_i$ is the fundamental group of
$\dom{\Omega_i}{\CS}$, the given isomorphism lifts to an
automorphism $\psi \co \CS \ra \CS$ such that $\psi \Omega_1
\psi^{-1} = \Omega_2$ in $\oper{Isom}(\CS)$, the group of
isometries of $\CS$.

By Cartan's theorem \cite[Chapter~IV]{Hel}
in characteristic zero and Tits' theorem \cite[Corollaries~5.9
and~5.10]{Tits} in positive characteristic, $\oper{Isom}(\CS) = $
the continuous automorphisms of $\PGL[d](F)$.

It is well known that $\Aut(\PGL[d](F))$ is a semidirect product
of the subgroup generated by (continuous) field automorphisms and
inner automorphisms,
by the subgroup of order $2$ generated by the `diagram
automorphism' $a \mapsto (\trans{a})^{-1}$. Replacing $D_2$ with
$\op{D_2}$ if necessary, we may ignore the last option. Therefore
$\psi$ extends to an automorphism of the matrix algebra
$\M[d](F)$.

\ver{4.2} Recall that $\Omega_i \sub G'_i(R_{0,T_i}) \sub G'_i(k)
\sub \PGL[d](F)$.
The quotient $\PGL[d](F)/\PSL[d](F)$ is abelian
torsion,
and $\Omega_i$ are finitely generated (as the fundamental groups
of compact objects). Therefore, we can replace the $\Omega_i$ by
finite index subgroups, which will be contained in $\PSL[d](F)$,
still satisfying the conjugation property above. Moreover, as
$\Aut(\PGL[d](F)) = \Aut(\PSL[d](F)) = \Aut(\SL[d](F)) = \Aut(\GL[d](F))$,
we can lift $\Omega_i$ to subgroups $\Omega_i^{(1)}$ of
$\SL[d](F)$, with $\psi \in \Aut(\GL[d](F))$ mapping
$\Omega_1^{(1)}$ to $\Omega_2^{(1)}$.

Let $H'_i$ denote the group of elements of norm one in
$\mul{D_i}$; in particular, $H'_i(F) = \SL[d](F)$. Then
$\Omega_i^{(1)} \sub H'_i(R_{0,T_i}) \sub H'_i(k) \sub \SL[d](F)$.
Let $N_i$ be the commensurator of $\Omega_i^{(1)}$ in $\GL[d](F)$.
We claim that $N_i$ is $\mul{D_i}(k)$ times the
scalar matrices.
\ver{4.3}
Indeed, an element $a \in N_i$ conjugates a finite index subgroup
$A$ of $\Omega_i^{(1)}$ into $\Omega_i^{(1)}$. Since
$\Omega_i^{(1)}$ is Zariski dense, so is $A$.  Therefore, the
$\bar{F}$-subalgebra of $\M[d](\bar{F})$ generated by $A$ is equal
to $\M[d](\bar{F})$ (where $\bar{F}$ is the algebraic closure).
This algebra is the scalar extension of the $k$-subalgebra
generated by $A$, which is therefore a $d^2$-dimensional
subalgebra of $D_i(k)$ (as $A \sub \mul{D_i}(k)$) and thus equal
to $D_i(k)$. Therefore, conjugation by $a$ is an automorphism of
$D_i(k)$. But by Skolem-Noether there is an element of
$\mul{D_i}(k)$ inducing the same automorphism, so up to scalar
matrices, $a \in \mul{D_i}(k)$.

Since $\psi$ maps $\Omega_1^{(1)}$ to $\Omega_2^{(1)}$, it maps
$N_1$ to $N_2$.  It also maps the commutator subgroup of $N_1$,
which is the commutator subgroup of $\mul{D_1}(k)$, to that of
$\mul{D_2}(k)$. It follows that $\psi$ maps the division subring
generated by commutators of $D_1(k)$ to the division subring
generated by
commutators of $D_2(k)$, but by
\Lref{subring} below, these are $D_1(k)$ and $D_2(k)$,
respectively. \ver{4.3}
Hence we proved that $\psi$ extends to an
isomorphism (of rings) from $D_1(k)$ to $D_2(k)$.

\ver{4.2} The restriction to $k$ of this isomorphism preserves
$\valF$ since it is induced by a continuous automorphism of $F =
k_{\valF}$. Also, it takes the ramification points $T_1$ of $D_1$
to the ramification points $T_2$ of $D_2$ (note that the opposite
algebra has the same set of ramification points).

\ver{4.2} We remark that in characteristic zero we could bypass
the argument with the commensurator, by observing that up to
complex conjugation, $\psi$ must be an algebraic map and hence
maps $G'_1(k)$ to $G'_2(k)$.

\ver{4.2} To show that $(3)\implies (2)$, we first assume $D_2(k)$
is isomorphic to $D_1(k)$. By tensoring with $F$, we obtain an
automorphism of $D_i(F) = \M[d](F)$, which is continuous by
assumption, and so induces an automorphism of $\CS$. Moreover it
maps the set $T_1$ to $T_2$, inducing an isomorphism from the
subring $R_{0,T_1}$ to $R_{0,T_2}$.  This extends to an
isomorphism from $G'_1(R_{0,T_1})$ to $G'_2(R_{0,T_2})$. Now, if
$\Omega_1 \leq G'_1(R_{0,T_1}) \leq \GL[d](F)$ is any finite index
torsion-free subgroup, then its image under the automorphism,
$\Omega_2$, is a finite index torsion-free subgroup of
$G'_2(R_{0,T_2})$, with $\dom{\Omega_1}{\CS}$ isomorphic to
$\dom{\Omega_2}{\CS}$.

If $D_2(k)$ is isomorphic to $\op{D_1(k)}$, then by tensoring with
$F$ we obtain an anti-automorphism of $\M[d](F)$, and the argument
follows throughout.
\end{proof}

\begin{Lem}\label{subring}
Let $D$ be a non-commutative division ring. Then the smallest
division subring of $D$ containing all multiplicative commutators
is $D$ itself.
\end{Lem}
\begin{proof}
First note that a division subring of $D$ containing all
commutators is invariant under conjugations (as $xyx^{-1} =
[x,y]y$). By the Cartan-Brauer-Hua
theorem \cite[p.~186]{Jac}
such a subring is either central or equal to $D$.

Assume all commutators in $D$ are central, and let $x \in D$ be a
non-central element.  Then the field generated by $x$ over
$\Cent(D)$ is invariant and therefore equal to $D$, so $D$ is
commutative, contrary to assumption.
\end{proof}

\ver{4.2} \Tref{noncommens} gives the condition for two manifolds
to be commensurable when they are defined from two division
$k$-algebras. To complete the proof of \Tref{isomanifold}, we need
to construct arbitrarily many division algebras with the same set
of ramification points, which are not isomorphic as rings.

\ver{4.3} The classification of division algebras over $k$ by
their invariants is crucial for that matter, so we add some (well
known) details on this topic (see
\cite[Chapsters~17--18]{Pierce}). Let $D$ be a division algebra,
with center $k$, a global field. Let $T$ denote the set of
valuations $\theta$ of $k$ such that $D \tensor k_{\theta}$ is not
a matrix algebra, where $k_{\theta}$ is the completion with
respect to $\theta$. Then $T$ is a finite set. Recall that by
Grunwald-Wang theorem \cite[Chapter~10]{AT}, there exists a cyclic
field extension $k_1$ of dimension $d$ over $k$, which is
unramified with respect to every valuation $\theta \in T$. {}From
Albert-Brauer-Hasse-Noether theorem it then follows that $D$
contains an isomorphic copy of $k_1$.

Let $\nu$ be a non-archimedean valuation of $k$. Let $\phi$ denote
the Frobenius automorphism of $k_1/k$ (\ie\ the generator of
$\Gal(k_1/k)$ which induces the Frobenius automorphism on the
residue field of $k_1$ over that of $k$). By Skolem-Noether
theorem, there is an element $z \in D$ such that conjugation by
$z$ induces $\phi$ on $k_1$. Then $z^d$ commutes with $k_1[z] =
D$, so $z^d$ is central, \ie\ $z^d \in k$. $D = k_1[z \subjectto
zaz^{-1} = \phi(a),\,z^d = b \in k]$.  If $k_2 \sub D$ is an
isomorphic copy of $k_1$, then again by Skolem-Noether there is an
element $w \in D$ such that $w k_1 w^{-1} = k_2$. Now $w z w^{-1}$
induces $\phi$ (or more precisely its isomorphic image) on $k_2$,
and the $d$-power of this element is $w z^d w^{-1} = z^d$. Of
course, $z$ is not unique; if $z'$ is another element inducing
$\phi$ on $k_1$, then $z'z^{-1}$ commutes with $k_1$.  So $z' =
az$ for some $a \in k_1$, since $k_1$ is its own centralizer.
Then, $(z')^d = \Norm[k_1/k](a)z^d$. However,
$\nu(\Norm[k_1/k](a))$ is divisible by $d$ since $k_1/k$ is
unramified, so $\nu(z^d) \pmod{d}$ is well defined. This number,
divided by $d$, is called the local invariant of $D$ at $\nu$
(\ver{4.2}the invariant is viewed as an element of $\Q/\Z$). The
completion $k_{\nu}$ is a splitting field of $D$ iff the invariant
is zero.

We remark that $D$ is determined (up to isomorphism as a
$k$-algebra) by its local invariants.  For every finitely
supported function $f \co \Vals\cup\set{\nu_0} \ra
\set{0,1/d,\dots,(d-1)/d}\subset \Q/\Z$ such that $\sum f(\nu)
\equiv 0 \pmod{1}$ and $f(\nu_0) = 0$, there is a central simple
$k$-algebra of degree $d$ which is matrices over $F = k_{\nu_0}$
and with $f(\nu)$  as the local invariant at $\nu$.

\begin{Prop}\label{invar}
Let $D$ and $D'$ be $k$-central division algebras which are
ramified at places $( \nu_1, \ldots, \nu_t)$ with invariants
$(a_1, \ldots, a_t)$ and $(\nu'_1, \ldots, \nu'_{t'})$ with
invariants $(a'_1, \ldots, a'_{t'})$, respectively.  Then $D$ and
$D'$ are isomorphic as rings if and only if $t = t'$ and there
exists a $k$-automorphism $\sigma$ such that (after reordering)
$\nu'_i = \nu_i \circ \sigma$, and $a_{i} = a'_i$.

Moreover, if $\bar{\sigma}$ is a fixed isomorphism $D \ra D'$
extending $\sigma$, then every isomorphism from $D$ to $D'$ is
equal to $\bar{\sigma}$ composed with an inner automorphism.
\end{Prop}
\begin{proof}
Assume $D'$ and $D$ are isomorphic as rings. The ring isomorphism
induces an isomorphism of the centers, which is an automorphism
$\s$ of $k$. Let $k_1$ denote an cyclic maximal subfield of $D$ as
before. If conjugation by $z$ induces $\phi$ on $k_1$, then
conjugation by its image $z'$ induces an automorphism $\phi'$ on
the image $k_1'$ of $k_1$, and of course $(z')^d = \s(z^d)$.
Therefore, the invariant of $D'$ at $\nu \circ \s^{-1}$ is
$(\nu\circ\s^{-1})((z')^d) = \nu(z^d)$, namely the invariant of
$D$ at the valuation $\nu$.  In the other direction, the given
$\sigma$ extends to an isomorphism of the rings translating the
ramification data of $D$ to that of $D'$.

Finally, if two isomorphisms $\psi,\psi' \co D \ra D'$ agree on
the center of $D$, then $\psi^{-1} \circ \psi'$ is an automorphism
of $D$ as a central division algebra.  By Skolem-Noether, it is an
inner automorphism.
\end{proof}

\begin{Cor}
\ver{4.3} Let $D_1$ and $D_2$ be $k$-division algebras of degree
$d$ with the same ramification points, such that for every
automorphism of $k$, the vector $a = (a_1,\dots,a_t)$ of ramified
invariants of $D_1$ is not permuted to the vector of ramified
invariants $b = (b_1,\dots,b_t)$ of $D_2$, nor to $-b$.

Let $\Omega_i \sub G'_i(R_{0,T})$ be finite index subgroups (where
$G'_i$ are defined as before). Then $\dom{\Omega_1}{\CS}$ and
$\dom{\Omega_2}{\CS}$ are not commensurable (where $\CS =
G(F)/K$). Equivalently, $\Omega_1$ and $\Omega_2$ are not
commensurable (after conjugation) as subgroups of
$\Aut(\PGL[d](F))$ (and in particular not as subgroups of
$\PGL[d](F)$).
\end{Cor}
\begin{proof}
Obviously the two claims are equivalent, for if $\Omega_0 =
\Omega_1 \cap g \Omega_2 g^{-1}$ is a finite index subgroup, then
$\dom{\Omega_0}{\CS} \isom \dom{g \Omega_0 g^{-1}}{\CS}$ is a
joint cover of finite index, and vice versa.

If such a joint cover exists, then $D_2$ must be isomorphic or
anti-isomorphic to $D_1$ as rings by \Tref{noncommens}, and this
is prohibited by the previous proposition.
\end{proof}

Before we continue, we first prove a Lemma.

\begin{Lem}
A global field $k$ has infinitely many orbits of valuations under
the action of its automorphism group.
\end{Lem}

\begin{proof}
Indeed, in zero characteristic there are infinitely many
valuations (covering the $p$-adic valuations for $p$ the rational
primes).  Since the Galois group over $\Q$ is finite, there are
finitely many valuations in each orbit.

In the case of a function field, $k$ is finite over $\F_q(t)$
where $\F_q$ is the maximal finite subfield. The degree of a
valuation is the dimension of its residue field over $\F_q$. If $p
\in \F_q[t]$ is a prime, then the degree of the $p$-adic valuation
on $\F_q(t)$ is equal to the degree of $p$ (as a polynomial), and
serves as a lower bound for the degree of an extension of this
valuation to $k$. Therefore, the degrees of valuations of $k$ are
unbounded. However, the automorphism group, when acting on the
valuations, does not change the degree.
\end{proof}

We can now finish the proof of Theorems \ref{isomanifold} and
\ref{isoComplex}.
Given
the number $m$,
let $t$ be an even integer such that $2^t/t \geq 2m$.
\ver{4.4} Choose a set $T \sub \Vals$ composed of $t$
non-archimedean valuations of $k$, which belong to different
orbits of the automorphism group of $k$. For every $e =
(e_1,\dots,e_t)$ with $e_i$ prime to $d$ which sum to zero modulo
$d$, let $D_e$ be the division
algebra over $k$ with
Brauer-Hasse-Noether invariants
$\frac{e_1}{d},\dots,\frac{e_t}{d}$ in the valuations of $T$,
which is unramified outside of $T$. Since $t$ is even we can
choose half of the invariants to be $1$ and half to be $-1$ ($1
\not \equiv -1$ as $d
>2$).  There are at least $\binom{t}{t/2} \geq 2^t/t \geq 2m$
options.
In each case, the sum of $e_i$ is zero modulo $d$. Furthermore we
set $e_1 = 1$, so there are at least $m$ vectors $e$ which are not
only different, but also not opposite to each other.

Since our valuations in $T$ are from different orbits and we have
chosen distinct invariants, there is no automorphism $\sigma$ of
$k$ which permutes the ramified valuations and their invariants.
This guarantees that the $D_e$ are not isomorphic to each other
(nor their opposites) as rings. Having chosen the division
algebras, take $U = U^{(r)}$ an open compact subgroup of
$G(\A_{\comp{T}-\set{\nu_0}})$ as in \Eq{Urdef}. Let $G'_e$ denote
the $k$-algebraic group $\mul{D_e}/\mul{Z}$ (where $Z$ is the
center), and take $\Gamma_{e} = \con{G'_e}{R_{0,T}}{I} = G'_e(k)
\cap G(F) U G'_e(\A_{T})$ (see \Eq{defGe}). We saw above that $X_e
= \dc{\Gamma_{e}}{G(F)}{K}$ are all isospectral, where $K$ is a
fixed
maximal compact subgroup of $G(F)$.
As a congruence subgroup, $\Gamma_e$ is of finite index in
$G'_e(R_{0,T})$, so the $X_e$ are non-commensurable by
\Tref{noncommens}.

\section{Proof of \Tref{basic}}\label{later}

Let $D$ be a division algebra of degree $d$ over the global field
$k$, and let $G'$ and $G$ denote the algebraic groups
$\mul{D}/\mul{Z}$ and $\PGL[d]$, respectively. Let $T$
denote the set of ramification points of $D$, and assume as before
that $D \tensor k_{\theta}$ is a division algebra for each $\theta
\in T$ (so the invariants of $D$ at $T$ are prime to $d$). Now let
$\theta \in T$. The
\defin{local \JLc} is a bijective correspondence, which maps
every irreducible, unitary representation $\rho'$ of
$G'(k_{\theta})$ to an irreducible, unitary square-integrable
representation $\rho = \JL_\theta(\rho')$ of $G(k_{\theta})$ (see
{\cite[p.~29]{HT}} or \cite{paperI} for details). Every such
representation of $G(k_{\theta})$ is either super-cuspidal (\ie\
the matrix coefficients are compactly supported modulo the
center), or is induced from a super-cuspidal representation of a
subgroup of smaller rank: $\rho = \Sp_s(\psi)$ where $\psi$ is a
super-cuspidal representation of $\PGL[d/s](k_{\theta})$, and
$\Sp_s$ is the construction of a `special representation'.  See
\cite{HT} or \cite[Subsection~2.5]{paperI}. We note that
$\Sp_1(\psi) = \psi$.

Let $\A$ denote the ring of \adeles\ over $k$. The automorphic
\rep{s} of $G(\A)$ are the irreducible \subrep{s} of
$\automrep[G]$ (as right $G(\A)$-modules), and likewise for
$G'(\A)$. Such a \rep{} decomposes as a tensor product of
\defin{local components} $\pi_{\nu}$ ($\pi=\otimes \pi_{\nu}$),
which are \rep{s} of $G(k_{\nu})$ (respectively $G'(k_{\nu})$) for
the various valuations $\nu$ of $k$. Moreover, $\pi$ is determined
by all but finitely many local components (this is `strong
multiplicity one' for $\PGL[d]$.)
.

The global \JLc\ maps an irreducible automorphic representation
$\pi' = \otimes \pi'_\nu$ of $G'(\A)$ to an irreducible
automorphic representation $\pi = \JL(\pi') = \otimes
\JL(\pi')_\nu$ of $G(\A)$ which occurs in the discrete spectrum
(see \cite[p.~195]{HT}).
If $\nu \not \in T$, then $\JL(\pi')_\nu \isom \pi'_\nu$.
For the valuations $\theta \in T$, the local correspondence takes
$(\pi')_{\theta}$ to $\JL_{\theta}((\pi')_{\theta}) = \Sp_s(\psi)$
for a suitable $s \divides d$ and $\psi$ super-cuspidal. Now in
the global map, the local component $\JL(\pi')_{\theta}$ is
isomorphic to either $\Sp_s(\psi)$ or to a certain representation
which we denote by
$\Cp_s(\psi)$, which is never square-integrable.
In both cases, $\psi$ and $s$ are uniquely determined by
$\pi'_\theta$.

By part 3 of \cite[Theorem~VI.1.1]{HT}, the image of the global
$\JL$ consists of the automorphic \rep{s} $\pi$ in the discrete
spectrum, such that for each $\theta \in T$, the local component
$\pi_{\theta}$ is either of the form $\Sp_s(\psi)$ or
$\Cp_s(\psi)$, a condition which depends on $T$ but is independent
of the invariants of $D$. Therefore the image of $\JL$ only
depends on the set $T$, and not on the invariants.
Let $\GAT = G'(\A_{\comp{T}})$ and $\Delta = G'(k) \sub \GAT$
(with the diagonal embedding). Recall that $\GAT =
G(\A_{\comp{T}})$, since $G'(k_{\nu}) = G(k_{\nu})$ for every $\nu
\not \in T$. First, we show that the only finite-dimensional
\subrep{}  of $\L2{\dom{\Delta}{\GAT}}$ is the trivial one. It is
well known that all finite-dimensional \subrep{s} are
one-dimensional, and so they come from the quotient
$\PGL[d](\A_{\comp{T}})/ \PSL[d](\A_{\comp{T}})$.  However, since
$T$ is not empty, $\PGL[d](k) \cdot \PSL[d](\A_{\comp{T}}) =
\PGL[d](\A_{\comp{T}})$, and so all such representations are
trivial.

Let us now define a map $\JL^{0}$, which takes an irreducible
infinite-dimensional \subrep{} $\pi''$ of
$\L2{\dom{\Delta}{\GAT}}$ to an irreducible infinite-dimensional
\subrep{} $\pi$ of $\L2{\dom{G(k)}{G(\A)}}$, which occurs in the
discrete spectrum. $\JL^{0}(\pi'')$ is defined as follows: Since
$\GAT = G'(\A)/G'(\A_{T})$ and
$\L2{\dc{G'(k)}{G'(\A)}{G'(\A_{T})}} \isom
\L2{\dom{G'(k)}{G'(\A)}}^{G'(\A_{T})} \sub
\L2{\dom{G'(k)}{G'(\A)}}$, $\pi''$ lifts to a \subrep{} $\pi'$ of
$\L2{\dom{G'(k)}{G'(\A)}}$ on which ${G'(\A_T)}$ acts trivially.
Namely, $(\pi')_{\theta}$ is trivial for every $\theta \in T$. By
the global \JLc, there is an automorphic \rep{} $\pi = \JL(\pi')$
of $G(\A)$ such that $(\pi')_{\nu} \isom \pi_{\nu}$ for every $\nu
\not \in T$. We set $\JL^{0}(\pi'') = \pi$. Now, for $\theta \in
T$, $\JL_{\theta}((\pi')_{\theta}) = \JL_{\theta}(1) =
Sp_d(\delta)$ is the Steinberg representation, for $\delta =
\absdot[]^{(1-d)/2}$. Moreover $\Cp_d(\delta)$ is the trivial
representation, and hence $\pi_{\theta}$ is either trivial or
Steinberg. But if $\pi$ is infinite-dimensional, $\pi_{\theta}$
cannot be trivial. Therefore, it is the Steinberg \rep{}.

The strong multiplicity one theorem
for $\mul{D}$ (see part 4 of \cite[Theorem~VI.1.1]{HT}) implies
that $\pi'$ has multiplicity one in $\L2{\dom{G'(k)}{G'(\A)}}$ and
hence $\pi''$ has multiplicity one in $\L2{\dom{\Delta}{\GAT}}$.
Moreover, it implies that $\JL^{0}(\pi'')$ determines $\pi''$. By
part 3 of \cite[Thmeorem~VI.1.1]{HT} and the discussion above, the
image of $\JL^{0}$ is composed of the discrete automorphic \rep{s}
$\pi$ of $\PGL[d](\A)$ such that for every $\theta \in T$, the
local component $\pi_{\theta}$ is a Steinberg \rep{.} In
particular, the image depends on $T$ but not on the specific
invariants of $D$.

\Tref{basic} follows immediately.  Since
$\L2{\dom{\Delta_1}{\GAT}}$ and $\L2{\dom{\Delta_2}{\GAT}}$ have
the same isomorphic images via the $\JL^{0}$ map in
$\L2{\dom{G(k)}{G(\A)}}$, they have isomorphic irreducible
\subrep{s}, each appearing with multiplicity one.

We should remark that \cite{HT} deals with
representations of $\GL[d]$ or $\mul{D}$ with a fixed central
character (and not $\PGL$ and $\mul{D}/\mul{Z}$), but this does
not change the argument as we can take the character to be
trivial.

\section{Generalization of \Tref{basic} and \Pref{basicprop}}

In this section we generalize \Tref{basic} and its main corollary.
The generalizations are not needed here, but we use them in a
subsequent paper on isospectral Cayley graphs \cite{paperIV}.

As before, let $k$ be a global field, $T$ a fixed finite set of
valuations, and $D_1,D_2$ central division algebras of degree $d$
over $k$, which are unramified outside of $T$, and remain division
algebras over the completions $k_{\theta}$ for $\theta \in T$. Let
$G'_i$ denote the multiplicative group of $D_i$ modulo center.
Recall that $G'_i(\A_{\comp{T}})$ can be identified with $\GAT =
G(\A_{\comp{T}})$. Let $J_i = G'_i(\A_T)$, so that $G'_i(\A) =
G'_i(\A_{\comp{T}}) J_i = \GAT J_i$.
Again $\Delta_i = G'_i(k)$ is embedded diagonally in $G'_i(\A)$.

\begin{Thm}\label{basic+}
Let $J_i^{0}$ be finite index subgroups of $J_i$ ($i = 1,2$), with
$\dimcol{J_1}{J_1^{0}} = \dimcol{J_2}{J_2^{0}}$.

Then $\L2{\dom{\Delta_1}{G'_1(\A)}}^{J_1^{0}}$ and
$\L2{\dom{\Delta_2}{G'_2(\A)}}^{J_2^{0}}$ (the spaces of
$J_i^{0}$-invariant functions on $\dom{\Delta_i}{G'_i(\A)}$) are
isomorphic as $\GAT$-\rep{s}.
\end{Thm}

Since $\L2{\dom{\Delta_i}{G'_i(\A)}}^{J_i} =
\L2{\dom{\Delta_i}{G'_i(\A_{\comp{T}})}}$, we obtain \Tref{basic}
as a special case by taking $J_i^{0} = J_i$. In fact,
\Tref{basic+} follows at once from \Tref{basic}: the natural
projection
$$\dom{\Delta_i}{G'_i(\A)} \ra \dom{\Delta_i}{G'_i(\A_{\comp{T}})}$$
has fiber $J_i$ over every point. Therefore
\begin{eqnarray*}
\L2{\dom{\Delta_i}{G'_i(\A)}}^{J_i^{0}} & = &
(\L2{\dom{\Delta_i}{G'_i(\A)}}^{J_i})^{[J_i:J_i^{0}]} \\
    & = &
    \L2{\dom{\Delta_i}{G'_i(\A_{\comp{T}})}}^{[J_i:J_i^{0}]}
    \end{eqnarray*}
where the first equality (for $\GAT$-modules) follows from the
fact that $J_i^{0}$ and $\GAT$ commute elementwise. As the indices
are equal by assumption, this completes the proof of
\Tref{basic+}.

Fix a valuation $\valF \not \in T$, and let $F = k_{\valF}$ be the
completion with respect to this valuation. Recall the rings
$R_{0,T}$ defined in \Eq{R0T}, and
\begin{equation}\label{R0}
R_{0} = \set{x \in k \suchthat \nu(x)\geq 0 \quad \mbox{for
every $\nu \neq \valF$}}.
\end{equation}
Of course, $R_0 \subset R_{0,T}$, and every ideal $I \normali
R_{0,T}$ induces the ideal $I \cap R_0$ of $R_0$. On the other
hand not every ideal of $R_0$ has this form (for example when $I
\normali R_0$ contains an element $a \in I$ with $\theta(a) = 0$
for some $\theta \in T$).

Fix embeddings $G'_i(k) \hookrightarrow \GL[t](k)$ (same $t$ for
both $i = 1,2$) such that over the completions $k_\theta$, $\theta
\in T$, $G'_i(k_{\theta}) = G'_i(\O_{\theta})$. This is possible
since $G'(k_{\theta})$ is compact for every $\theta$ in the finite
set $T$. Set $G'_i(R_0) = G'_i(k) \cap \GL[t](R_0)$, and define
the congruence subgroups $\con{G'_i}{R_0}{I}$ as in \Eq{condef}.
We obtain the following extension of \Pref{basicprop}.

\begin{Prop}\label{basicprop+}
Let $0 \neq I \normali R_{0,T}$ be an ideal, and let $J_i^{0}$ be
finite index subgroups of $J_i = G'_i(A_T)$ with
$\dimcol{J_1}{J_1^0} = \dimcol{J_2}{J_2^{0}}$. Let $\Lambda_i =
G'(k) \cap G(F) U^{(r)} J_i^{0}$ be the corresponding discrete
subgroup, where $U^{(r)}$ is defined as in \eq{Urdef}. Then
\begin{equation*}
\L2{\dom{\Lambda_1}{G(F)}} \isom \L2{\dom{\Lambda_2}{G(F)}}
\end{equation*}
as representations of $G(F) = \PGL[d](F)$.
\end{Prop}
\begin{proof}
Let $\Delta_i = G'_i(k)$ as before. \Tref{basic+} implies that
$$\L2{\dom{\Delta_1}{G'_1(\A)}}^{U^{(r)} J_1^{0}} \isom \L2{\dom{\Delta_2}{G'_2(\A)}}^{U^{(r)}
J_2^{0}},$$ so
$$\L2{\dc{\Delta_1}{G'_1(\A)}{U^{(r)}J_1^0}} \isom
\L2{\dc{\Delta_2}{G'_2(\A)}{U^{(r)}J_2^0}}.$$

On the other hand, combining \Lref{finqd} and the fact that
$\dimcol{J_i}{J_i^{0}}$ is independent of $i$, we see that
$$G'_i(\A) / \Delta_i G(F) U^{(r)} J_i^{0}$$
is independent of $i$. As in the proof of \Pref{basicprop}, this
implies that
$$\L2{\dc{\Delta_i}{\Delta_i G(F) U^{(r)} J_i^{0}}{U^{(r)}J_i^0}}$$
is independent of $i$. But this space is isomorphic to
$\L2{\dom{\Lambda_i}{G(F)}}$ by the choice of $\Lambda_i$.
\end{proof}

This proposition generalizes \Pref{basicprop}, providing
isospectrality for principal congruence subgroups with respect to
ideals of $R_0$, rather than only ideals of $R_{0,T}$. (However it
also covers non-principal congruence subgroups).
\begin{Cor}
Let $0 \neq I \normali R_{0}$. There is an isomorphism
\begin{equation*}
\L2{\dom{\con{G'_1}{R_{0}}{I}}{G(F)}} \isom
\L2{\dom{\con{G'_2}{R_{0}}{I}}{G(F)}},
\end{equation*}
as representations of $G(F) = \PGL[d](F)$.
\end{Cor}
\begin{proof}
Define a function $r \co \Vals \minusset \set{\valF} \ra \N \cup
\set{0}$ as in \Eq{rdef}, and let $U^{(r)}$ be the group defined
in \Eq{Urdef}. Take $J_i^{0} = \prod_{\theta \in T}
\con{G'_i}{\O_{\theta}}{P_{\theta}^{r_{\theta}}}$, so that
$\con{G'_i}{R_0}{I} = G'_i(k) \cap G(F) U^{(r)} J_i^0$. {}From the
detailed description of the quotients
$G'(k_{\theta})/\con{G'_i}{\O_{\theta}}{P_{\theta}^{r_{\theta}}}$
in \cite[Proposition~1.5]{PrasadRagh}, one sees that the index of
$J_i^{0}$ in $G'_i(\A_T)$ does not depend on $i$, so the
proposition applies.
\end{proof}

\def\US{A.~Lubotzky, B.~Samuels and U.~Vishne}
\def\BIBPapI{\US, {\it Ramanujan complexes of type $\tilde{A_d}$},
Israel J. of Math., to appear.}
\def\BIBPapII{\US, {\it Explicit constructions of Ramanujan complexes},
European J. of Combinatorics, to appear.}
\def\BIBPapIII{\US, {\it Division algebras and non-commensurable isospectral manifolds},
preprint.}
\def\BIBPapIV{\US, {\it Isospectral Cayley graphs of some simple groups},
preprint.}

\end{document}